\documentclass[12pt]{article}

\usepackage[margin=1in]{geometry} 
\usepackage{graphicx}             
\usepackage{amsmath}               
\usepackage{amsfonts}             
\usepackage{amsthm}               
\usepackage{amssymb,latexsym,hyperref}

\newtheorem{thm}{Theorem}[section]
\newtheorem{lem}[thm]{Lemma}
\newtheorem{claim}[thm]{Claim}
\newtheorem{defn}[thm]{Definition}

\def\tree{\mathcal{T}}
\def\env{\text{env}}

\begin{document}

\nocite{*}

\title{Avoiding Rational Distances}

\author{Ashutosh Kumar\footnote{University of Wisconsin Madison, Wisconsin 53706; email: akumar@math.wisc.edu}}

\maketitle

\begin{abstract}
 
We show that for any set of reals $X$ there is a $Y \subseteq X$ such $X$ and $Y$ have same Lebesgue outer measure and the
distance between any two distinct points in $Y$ is irrational.

\end{abstract}

\section{Introduction}

P\'{e}ter Komj\'{a}th has asked the following question in \cite{kom1}: Let $X$ be a subset of Euclidean space $\mathbb{R}^n$.
Is there always a $Y \subseteq X$ such that $X$ and $Y$ have same outer measure and the distance between any two
distinct points of $Y$ is irrational? In \cite{kom2} he showed that $\mathbb{R}^n$ can be colored by countably many colors
such that the distance between any two points of the same color is irrational. It follows that one can always find a subset
of positive outer measure that avoids rational distances. Under the assumption that there is no weakly inaccessible cardinal
below the continuum, he also showed in \cite{kom1} that in dimension one we can always find a subset $Y$ of full outer 
measure in $X$, avoiding rational distances. Moti Gitik and Saharon Shelah showed the following in \cite{GS1}, \cite{GS2}: For any sequence 
$\langle A_n : n \in \omega \rangle$ of sets of reals, there is disjoint refinement of full outer measure; i.e., there is a sequence
$\langle B_n : n \in \omega \rangle$ of pairwise disjoint sets such that $B_n \subseteq A_n$ and they have the same outer measure.
It follows that one can omit integer distances in dimension one while preserving outer measure. Their proof employs one of their results
about forcing with ideals that says that forcing with a sigma ideal cannot be isomorphic to a product of Cohen and Random forcings. Here
we answer Komj\'{a}th's question positively in dimension one.

\textbf{Acknowledgements}:  I am grateful to my advisor Arnold Miller for several useful discussions. 
I would also like to thank P\'{e}ter Komj\'{a}th for informing me about this question.

\section{A theorem of Gitik and Shelah}

Suppose $A$ is a subset of $\mathbb{R}^n$. We say that $B \subseteq A$ is full in $A$ if $\env(A) = \env(B)$
where by $\env(X)$ we denote a $G_{\delta}$ measurable envelope of $X$; i.e., $\env(X)$ is a $G_{\delta}$ set containing $X$ such that
the inner measure of $\env(X) \backslash X$ is zero. If the outer measure of $A$ is finite this is equivalent to saying
that $A$ and $B$ have same outer measure. 

Let $\tree$ be a subtree of $\omega^{< \omega}$ such that every node in 
$\tree$ has at least two children; i.e., for every $\sigma \in \tree$, $|\{n \in \omega : \sigma n \in \tree\}| \geq 2$.

\begin{defn} 

Call a family $\langle A_{\sigma} : \sigma \in \tree \rangle$
of subsets of a set $A$, a full tree on $A$ if:

\begin{itemize}

\item $A = A_{\phi}$, and for every $\sigma \in \tree$,

\item $A_{\sigma}$ is a disjoint union of $A_{\sigma n}$'s where $\sigma n \in \tree$

\item $A_{\sigma}$ is full in $A$.

\end{itemize}
 
\end{defn}

The following application of Theorem 2.3 is implicit in \cite{GS1}:

\begin{thm} 

Let $A \subseteq \mathbb{R}^n$ and let $\langle A_{\sigma} : \sigma \in \tree \rangle$
be a full tree on $A$. Then there is a $B \subseteq A$ full in $A$ such that for every $\sigma \in \tree$, 
$A_{\sigma} \backslash B$ is full in $A_{\sigma}$.

\end{thm}

This theorem is a consequence of the following theorem in \cite{GS2}:

\begin{thm} 

Suppose $I$ is a sigma ideal over a set $X$. Then forcing with $I$ cannot be isomorphic to Cohen $\times$ Random.

\end{thm}

Let us explain how Theorem 2.2 follows from Theorem 2.3. It is clearly enough to show that there is a non null $B \subseteq A$ 
such that $A_{\sigma} \backslash B$ is full in $A_{\sigma}$ for every $\sigma \in \tree$, for then we can subtract 
$\env(B)$ from every node of our tree and repeat until we exhaust $\env(A)$. Suppose that this fails so that for every 
non null $B \subseteq A$, there is some $\sigma \in \tree$ such that $\env(A_{\sigma})$ is strictly larger than 
$\env(A_{\sigma} \backslash B)$. Consider the map that sends every positive outer measure subset $B \subseteq A$ to the supremum,
in the complete Boolean algebra Cohen $\times$ Random, of all pairs $(\sigma, E)$ where $\sigma \in \tree$ and $E$ is a 
positive measure Borel subset of $\env(A)$ such that $E$ is disjoint with $\env(A_{\sigma} \backslash B)$. This gives a dense 
embedding from $\mathcal{P}(A)\slash \text{Null}$ to Cohen $\times$ Random contradicting the fact that they cannot be forcing
isomorphic.

\section{The main result}

\begin{thm} 

Let $X \subset \mathbb{R}$ be a set of positive outer measure. Then there is a
$Y \subseteq X$ such that $Y$ is full in $X$ and the distance of any pair of distinct points
in $Y$ is irrational.

\end{thm}

Proof of Theorem 3.1: Let $|X| = \kappa$. Let $X_0 = \langle x_{\alpha} : \alpha < \kappa \rangle$ be a set of representatives from the
partition on $X$ induced by the relation $x \sim y$ iff $x - y$ is rational. Let $\langle r_n : n \geq 1 \rangle$ be a list
of all nonzero rationals. For each $n \geq 1$, let $f_{n}: X_0 \rightarrow \mathbb{R}$ be
defined by $f_n(x_{\alpha}) = x_{\alpha} + r_n$, if $x_{\alpha} + r_n \in X$, otherwise $f_n(x_{\alpha}) = x_{\alpha}$,
also put $X_n = \text{range}(f_n)$. Let $f_0$ be identity on $X_0$. For $n > m \geq 1$, 
let $F^{m}_{n} = f_n \circ f_m^{-1}: X_m \rightarrow X_n$.

We will inductively define a sequence $\langle K_n : n \geq 0 \rangle$ of pairwise disjoint subsets of
$\kappa$ such that $X_n \upharpoonright K_n = \langle f_n(x_{\alpha}) : \alpha \in K_n \rangle$ is full in $X_n$.
Theorem 3.1 will immediately follow. We'll need the following lemma:

\begin{lem} 

Let $n > m \geq 0$, and $Y$ be a positive outer measure subset of $X_m$.
Then there is a partition $\{Y_i : 1 \leq i \leq k\}$ of $Y$, such that for every $i$,

\begin{itemize}

\item $Y_i$ is full in $Y$ and

\item for every $W \subseteq Y_i$, if $W$ is full in $Y_i$, then 
$F^{m}_{n}[W]$ is full in $F^{m}_{n}[Y_i]$. In this case we say that $F^{m}_{n} \upharpoonright Y_i$
is fullness preserving.

\end{itemize}

\end{lem}

Proof of Lemma 3.2: We will make several uses of the following result of Luzin: Any set of reals $X$ can be partitioned into two full subsets (\cite{luz}).
Note that $F = F^{m}_{n} \upharpoonright Y$ translates a finite list of pieces of $Y$. Let
$Y = T_1 \sqcup T_2 \sqcup \dots \sqcup T_k$ where $T_i$ is translated by some rational $r_i$ (possibly, some $r_i = 0$).
Use induction on $k$. If $k = 1$ $Y_0 = Y$ works. So assume $k = l + 1$. Let $Z = \bigcup \{T_i : 1 \leq i \leq l\}$.
Let $\{Z_i : 1 \leq i \leq l\}$ be a partition of $Z$ such that each $Z_i$ is full in $Z$ and $F \upharpoonright Z_i$
is fullness preserving. Let $E_1 = \env(Z)$, $E_2 = \env(Y_k)$ and $D = E_1 \bigcap E_2$. Let $W_1$, $W_2$ be
a partition of $Z_1 \bigcap (E_1 \backslash D)$ into two full subsets. Let $\{V_j : 1 \leq j \leq k\}$ be a partition
of $Y_k \bigcap (E_2 \backslash D)$ into $k$ full subsets. Set $Y_1 = W_1 \bigcup (Z_1 \bigcap D) \bigcup V_1$.
For $2 \leq i \leq l$, put $Y_i = Z_i \bigcup V_i$ and let $Y_k = W_2 \bigcup (D \bigcap Y_k) \bigcup V_k$. Then $\{Y_i : 1 \leq i \leq k\}$
is a partition of $Y$ with the required properties.

\begin{claim} 

There exists $K_0 \subseteq \kappa$, such that $X_0 \upharpoonright K_0 = \{x_{\alpha} : \alpha \in K_0 \}$ is full in
$X_0$ and for every $n \geq 1$, $X_n \upharpoonright (\kappa \backslash K_0)$ is full in $X_n$.

\end{claim}

Proof of Claim 3.3: Using Lemma 3.2, construct a full tree $\langle Y_{\sigma} : \sigma \in 2^{< \omega} \rangle$ on $Y = X_0$ such that

\begin{itemize}
 
\item for each $\sigma \in 2^{n}$, $n \geq 1$ and for each $1 \leq i \leq n$, $f_{i} \upharpoonright Y_{\sigma}$
is fullness preserving

\end{itemize}

Now Theorem 2.2 will imply that there is some $K_0 \subseteq \kappa$ such that $X_0 \upharpoonright K_0$ is full in 
$X_0$ and for every $\sigma \in 2^{< \omega}$, $Y_{\sigma} \upharpoonright (\kappa \backslash K_0)$ is full in
$Y_{\sigma}$. Fix any $n \geq 1$ and note that for each $\sigma \in 2^n$, $f_{n} \upharpoonright Y_{\sigma}$ is fullness
preserving so that $f_{n}[Y_{\sigma} \upharpoonright (\kappa \backslash K_0)]$ is full in $f_{n}[Y_{\sigma}]$.
It follows that $X_n \upharpoonright \kappa \backslash K_0 = \bigcup \{f_{n}[Y_{\sigma} \upharpoonright 
(\kappa \backslash K_0)] : \sigma \in 2^n \}$ is full in $\bigcup \{f_{n}[Y_{\sigma}] : \sigma \in 2^n\} = X_n$.

Now suppose we have already obtained pairwise disjoint subset $\{K_i : 0 \leq i \leq n\}$ of subsets of $\kappa$ such that

\begin{itemize}
 
\item for each $0 \leq i \leq n$, $X_i \upharpoonright K_i$ is full in $X_i$

\item for each $j > n$, $X_j \upharpoonright (\kappa \backslash \bigcup \{K_i : 1 \leq i \leq n\})$ is full in $X_j$.

\end{itemize}

Then we can (using Lemma 3.2 and Theorem 2.2) obtain $K_{n+1}$ as in the proof of Claim 3.3 above.

\end{document}